\documentclass[12pt]{amsart} 
\usepackage{amsmath}
\usepackage{amsthm}
\usepackage{amssymb}
\usepackage{amscd}
\usepackage[dvips]{graphics} 
\usepackage{epsfig}
\usepackage{graphicx}
\usepackage{pgf}
\usepackage{epstopdf}

\newtheorem{theorem}{Theorem}[section]
\newtheorem{lemma}[theorem]{Lemma}
\newtheorem{proposition}[theorem]{Proposition}

\newtheorem{question}{Question}

\newtheorem*{theoremC}{Theorem~\ref{thmain}}

\theoremstyle{definition}

\newtheorem{remark}[theorem]{Remark}

\begin{document}

\def \k {K_{(a, b)}}
\def \l {L(p, q)}
\def \kp {K_{(a', b')}}

\def \L {L_{(a, b)}}
\def \Lpp {L_{(a', b')}}
\def \Lp {L'_{(a, b)}}

 \title{Legendrian Torus Knots in Lens Spaces}

 \author{S\.{i}nem Onaran}

 \address{Hacettepe University Department of Mathematics 06800 Beytepe-Ankara Turkey}
  \email{sonaran@hacettepe.edu.tr}
   \subjclass{57R17}
\keywords{Legendrian knots, contact structures, lens spaces}

\begin{abstract}
In this note, we first classify all topological torus knots lying on the Heegaard torus in lens spaces, and then we study Legendrian representatives of these knots. We classify oriented positive Legendrian torus knots in the universally tight contact structures on the lens spaces up to contactomorphism. 
\end{abstract}
\maketitle
 \setcounter{section}{0}

\section{Introduction} 

A \emph{Legendrian knot} in a contact $3$-manifold is a knot which is everywhere tangent to the contact planes. Legendrian knots are natural objects in contact $3$-manifolds and they are used to distinguish contact structures \cite{Kanda1}, to detect topological properties of knots \cite{R} and to detect overtwistedness of contact structures \cite{EH2}. 

\par There has been some recent progress in the classification of Legendrian knots in the standard tight contact structure on $S^3$ after the classification of Legendrian unknots made by Eliashberg and Fraser \cite{EFr} and the classification of Legendrian torus knots and the figure eight knot made by Etnyre and Honda \cite{EH}. Legendrian knots in a cabled knot type in $S^3$ are studied in \cite{EH3} and complete classification is given in \cite{T}. Recently, Legendrian twists knots are classified in \cite{ENV}. Legendrian knots in $3$--manifolds other than the $3$--sphere $S^3$ are also studied. For example, in \cite{Pao}, Legendrian linear curves in the $3$--torus $T^3$ are classified and in \cite{EB, GO}, Legendrian rational unknots in lens spaces are classified. 

\par In this note, we employ the classification scheme of Etnyre and Honda in \cite{EH} to Legendrian knots that are rationally null-homologous. We focus on Legendrian torus knots in lens spaces. For relatively prime integers $a$, $b$, an oriented simple closed curve that wraps $a$ times in the meridional direction and $b$ times in the longitudinal direction on the Heegaard torus is called an \emph{$(a, b)$--torus knot} in the lens space. We call $(a, b)$--torus knots with $a, b> 0$ \emph{positive} torus knots. In Section $2$, we study the topological properties of all torus knots in lens spaces. We find a constraint on when a torus knot is null-homologous. We compute the group of a torus knot which is defined as the fundamental group of its complement. By studying the diffeotopy group of lens spaces, we completely classify all torus knots up to isotopy. Lastly, we construct a rational Seifert surface for a torus knot and we calculate its Euler characteristic.  

In section $3$, we give a review of the basic concepts in convex surface theory and we fix the notation. In Section $4$, by using convex surface theory tools, we study Legendrian representatives of positive torus knots in the universally tight contact structures on the lens spaces.  We define the rational Legendrian knot invariants, \cite{F, BG, EB, GO}. We compute the rational Thurston-Bennequin invariants and the rational rotation numbers of oriented positive Legendrian torus knots by using the rational Seifert surface we constructed for torus knots. By following the strategy outlined in \cite{EH}, we first classify oriented positive Legendrian torus knots with maximal rational Thurston-Bennequin invariant. Then, we show that all oriented positive Legendrian torus knots with non-maximal rational Thurston-Bennequin invariant destabilize. We prove:

\begin{theoremC} Two oriented positive Legendrian torus knots in a universally tight contact structure on a lens space are uniquely realized up to contactomorphism if and only if their oriented knot types, rational Thurston-Bennequin invariants and rational rotation numbers agree.
\end{theoremC}

Legendrian knots may be classified up to contact isotopy or up to global contactomorphism. By the work of Eliashberg we know that the group of co-orientation preserving contactomorphisms of the standard tight $S^3$ is connected, {\cite[Theorem 2.4.2]{Elias}}. Therefore, for Legendrian knots in the standard tight $S^3$, these two classifications are equivalent. However, for arbitrary tight contact closed $3$--manifolds the group of co-orientation preserving contactomorphisms is not well understood. In particular, nothing is known for tight contact lens spaces in general.
\begin{question} Is the group of co-orientation preserving contactomorphisms of universally tight contact structures on lens spaces connected?
\label{Q1}
\end{question}
We want to remark that a positive answer to Question~\ref{Q1} together with Theorem~\ref{thmain} provides us the classification of Legendrian torus knots up to Legendrian isotopy in the universally tight contact structures on the lens spaces.

A \emph{transverse knot} in a contact 3-manifold is a knot which is everywhere
transverse to the contact planes. There are two types of classical invariants
for null-homologous transverse knots; the knot type and the self-linking number. The self-linking number may be generalized for rationally null-homologous transverse knots,
\cite{EB}. By \cite[Theorem 2.10]{EH}, two transverse knots in a contact $3$-manifold are transversely isotopic if and only if their Legendrian push offs are Legendrian isotopic after
each has been negatively stabilized some number of times. As a consequence of Theorem ~\ref{thmain} we have
\begin{theorem} Two positive transverse torus knots in a universally tight contact structure on a lens space are uniquely realized up to contactomorphism if and only if their knot types and rational self-linking numbers agree. 
\end{theorem}

\section{Topological torus knots in lens spaces} 

\par For fixed relatively prime integers $p > q > 0$, let $(V_1, V_2)$ be the genus $1$ Heegaard splitting of a lens space $L(p, q)$ which is described as 
\begin{equation*}
L(p, q) = V_1 \cup _{\varphi} V_2
\end{equation*}
where $V_1$ and $V_2$ are both $D^2\times S^1$. Let $\mu_i$ and $\lambda_i$ be a meridian and longitude pair for $V_i$, $i = 1,2$. The gluing map $\varphi : \partial V_1 \rightarrow \partial V_2$ is an orientation reversing map given in standard longitude-meridian coordinates on the torus by the matrix 

\begin{displaymath} \left(\begin{array}{ccc} -q&q'\\p&p' \end{array}\right) \end{displaymath} 
with $pq' + qp' = 1$, $p', q' \in \mathbb{Z}$. In particular, the image of the meridian $\mu_1$ of $\partial V_1$ is the curve $-q\mu_2 + p\lambda_2$ in $\partial V_2$. Note that $H_1(\l, \mathbb{Z}) \cong \mathbb{Z}/ p\mathbb{Z}$ is generated by $[\lambda_2]$. Therefore, any oriented knot $K$ in $\l$ represents $b[\lambda_2]$ for some integer $b$. Any knot in a lens space $\l$ is rationally null-homologous and has an order. The order $r$ of $K$ is defined to be the order of $[K]$, and hence $r = p/gcd(p,b)$. 

\par By Theorem 1 of \cite{Bon}, one can fix the Heegaard torus up to isotopy in a lens space. This allows us to define torus knots on the Heegaard torus $\partial V_2$. For relatively prime integers $a$, $b$, an oriented simple closed curve $\k$ that wraps $a$ times in the meridional direction and $b$ times in the longitudinal direction on $\partial V_2$ is called an \emph{$(a, b)$--torus knot} in the lens space $\l$. For a knot $\k$ of order $r$ in $\l$, $p \mid rb$ and $\k$ is null-homologous if and only if $p \mid b$. 

\begin{proposition} Let $\k$ be an $(a, b)$--torus knot on the Heegaard torus $\partial V_2$ in lens space $\l$.
\begin{enumerate} 
\item The group of a torus knot $\k$ can be presented as 
\begin{equation*}
\pi_1{(\l - \k)} = < u, v \mid u^b = v^{pa+qb} >.
\end{equation*}
\item Two torus knots $\k$ and $\kp$ have isomorphic groups if and only if $|b| = |b'|$ and $|pa + qb| = |pa' + qb'|$ or $|b| = |pa' + qb'|$ and $|b'| = |pa + qb|$.
\end{enumerate}
  \label{prop1}
\end{proposition}
\begin{proof} The complement of a neighbourhood $\nu(\k)$ of the torus knot $\k$ in $\l$ is the union of two solid tori glued along an annulus $A$ where the core $C$ of the annulus $A$ is isotopic to the torus knot $\k$. Namely, $\l \setminus  \nu(\k) = \tilde{V}_1 \cup \tilde{V}_2$ where $\tilde{V}_i = \overline{V_i \setminus \nu(\k)}$, $i = 1,2$, are two solid tori glued along the annulus $A = (\l \setminus \nu(\k)) \cap \partial V_2$. 
\par Let $\tilde{\mu}_i$ and $\tilde{\lambda}_i$ be a meridian and longitude pair for $\tilde{V}_i$ where $\tilde{\mu}_i$ and $\tilde{\lambda}_i$ represent the trivial element and a generator of $\pi_1(\tilde{V}_i)$, respectively. Note that the homotopy class $[C] = {[\tilde{\mu}_1]}^{- p'a + q'b}{[\tilde{\lambda}_1]}^{pa+qb} = {[\tilde{\lambda}_1]}^{pa+qb}$ since $\k$ is on $\partial V_2$ and 
\begin{displaymath} \left(\begin{array}{ccc} -q&q'\\p&p' \end{array}\right)^{-1}\left(\begin{array}{ccc} a\\b \end{array}\right) =   \left(\begin{array}{ccc} -p'&q'\\p&q \end{array}\right)\left(\begin{array}{ccc} a\\b \end{array}\right) = \left(\begin{array}{ccc} - p'a + q'b\\pa +qb \end{array}\right) \end{displaymath} 
for $pq' + qp' = 1$, $p', q' \in \mathbb{Z}$. Also, $[C] = {[\tilde{\mu}_2]}^a{[\tilde{\lambda}_2]}^b = {[\tilde{\lambda}_2]}^b$. Then, by Seifert-van Kampen theorem,  
\begin{equation*}
\pi_1{(\k)} = < u, v \mid u^b = v^{pa+qb} >
\end{equation*}
 where $u = [\tilde{\lambda}_1]$ and $v = [\tilde{\lambda}_2]$. This proves ($1$).

\par The subgroup $<u^b>$ is the center of the knot group $\pi_1{(\k)}$ and $\pi_1{(\k)}/<u^b> = \mathbb{Z}_{|b|} \ast \mathbb{Z}_{|pa+qb|}$. Note that $u$ and $v$ generate non-conjugate maximal finite cyclic subgroups of order $|b|$ and $|pa+qb|$ of $\mathbb{Z}_{|b|} \ast \mathbb{Z}_{|pa+qb|}$, respectively. Therefore, if $\k$ and $\kp$ have isomorphic groups, then $|b| = |b'|$ and $|pa + qb| = |pa' + qb'|$ or $|b| = |pa' + qb'|$ and $|b'| = |pa + qb|$. Now if $|b| = |b'|$ and $|pa + qb| = |pa' + qb'|$ or $|b| = |pa' + qb'|$ and $|b'| = |pa + qb|$, then from ($1$) it is straightforward that $\k$ and $\kp$ have isomorphic groups, proving ($2$).
\end{proof}

\begin{lemma}\begin{enumerate}
\item If $\k$ and $\kp$ are two null-homologous oriented torus knots in $\l$ that have isomorphic groups, then $(a', b')$ is equal to one of the following pairs: 
\begin{itemize}
\item $A = (a, b)$, $-A = (-a, -b)$, \\ $B = (\frac{-2qb-pa}{p}, b)$, $-B = (\frac{2qb+pa}{p}, -b)$,\\$C = (\frac{b-qpa-q^2b}{p}, pa+qb)$, $-C = (\frac{-b+qpa+q^2b}{p}, -pa-qb)$,\\$D = (\frac{b+qpa+q^2b}{p}, -pa-qb)$, $-D = (\frac{-b-qpa-q^2b}{p}, pa+qb)$.
\end{itemize}
\medskip
\item If $\k$ and $\kp$ are two rationally null-homologous but not null-homologous oriented torus knots in $\l$ that have isomorphic groups, then $(a', b')$ is equal to one of the following pairs in the following cases: 
\begin{itemize}
\item $A = (a, b)$, $-A = (-a, -b)$ if $p \neq 2$ and $q^2 \not\equiv \pm 1$ (mod $p$),
\item $A = (a, b)$, $-A = (-a, -b)$, $C = (\frac{b-qpa-q^2b}{p}, pa+qb)$, $-C = (\frac{-b+qpa+q^2b}{p}, -pa-qb)$ if $p \neq 2$ and $q^2 \equiv 1$ (mod $p$),
\item $A = (a, b)$, $-A = (-a, -b)$, $D = (\frac{b+qpa+q^2b}{p}, -pa-qb)$, $-D = (\frac{-b-qpa-q^2b}{p}, pa+qb)$ if $p \neq 2$ and $q^2 \equiv -1$ (mod $p$),
\item $A = (a, b)$,  $ -A = (-a, -b)$, $B = (-b-a, b)$, $-B = (b+a, -b)$, $C = (-a, 2a+b)$, $-C = (a, -2a-b)$, $D = (a+b, -2a-b)$, $-D = (-a-b, 2a+b)$ if $p = 2$.
\end{itemize}
\end{enumerate}
  \label{lemma1}
\end{lemma}
\begin{proof} By Proposition ~ \ref{prop1}$(2)$ we know that $\k$ and $\kp$ have isomorphic groups if and only if $|b| = |b'|$ and $|pa + qb| = |pa' + qb'|$ or $|b| = |pa' + qb'|$ and $|b'| = |pa + qb|$. Case $(1)$ follows from the analysis of these cases using that $p \mid b$. For Case $(2)$, we know that if $\k$ is not null-homologous then $p \nmid b$. Therefore, when $p \neq 2$ the cases $(\frac{-2qb-pa}{p}, b)$, $(\frac{2qb+pa}{p}, -b)$ do not occur and the cases $(\frac{b-qpa-q^2b}{p}, pa+qb)$, $(\frac{-b+qpa+q^2b}{p}, -pa-qb)$ occur only if $p \mid (1 - q^2 )$. Similarly, the cases $(\frac{b+qpa+q^2b}{p}, -pa-qb)$, $(\frac{-b-qpa-q^2b}{p}, pa+qb)$ occur only if $p \mid (1+q^2)$. The case when $p = 2$ and hence $q = 1$ is clear.
\end{proof}

\par Let us now classify all topological torus knots in lens spaces up to isotopy. Recall that two knots $K_1$ and $K_2$ in a $3$--manifold $M$ are isotopic if there is a diffeomorphism $g : M \rightarrow M$ such that $g(K_1) = K_2$ and $g$ is isotopic to the identity map.

\begin{theorem} The oriented torus knot $\k$ is isotopic to $\kp$ in $\l$ if and only if $(a', b')$ is an element of one of the following: 
\begin{enumerate}
\item $\{(a, b)\}$ if $q \neq 1$ or $p-1$,
\item $\{(a, b), (\frac{b-qpa-q^2b}{p}, pa+qb)\}$ if $p\neq 2$ and $q = 1$ or $p-1$,
\item $\{ (a, b), (-a, -b), (-a, 2a+b), (a, -2a-b) \}$ if $p = 2$.
\end{enumerate}
\label{class1}
\end{theorem}
\par For the proof of Theorem ~ \ref{class1} we need the following theorem:

\begin{theorem}[{\cite[Theorem 3]{Bon}}] The group of isotopy classes of diffeomorphisms of $\l$ for $p \geq 2$ is given by
\begin{enumerate}
\item $\mathbb{Z}_2$ with generator $\tau$ if $q^2 \not\equiv \pm 1$ (mod $p$), 
\item $\mathbb{Z}_2 \oplus \mathbb{Z}_2$ with generator $\tau$ and $\sigma_{+}$ if $q^2 \equiv 1$  and $q \not\equiv \pm 1$ (mod $p$),
\item $\mathbb{Z}_2$ with generator $\tau$ if $q \equiv \pm 1$ (mod $p$) and $p\neq 2$,
\item $\mathbb{Z}_4$ with generator $\sigma_{-}$ if $q^2 \equiv -1$ (mod $p$) and $p \neq 2$,
\item $\mathbb{Z}_2$ with generator $\sigma_{-}$ if $p = 2$.
\end{enumerate}
\label{thmlens}
\end{theorem}
Let $(V_1, V_2)$ be the genus $1$ Heegaard splitting of the lens space $L(p, q)$ defined as above. In Theorem ~\ref{thmlens}, the diffeomorphism $\tau$ preserves each of the solid tori $V_i = D^2 \times S^1$ and acts by a complex conjugation on each factor (as viewed in $\mathbb{C}$) of each Heegaard torus. Note that $\tau$ always exists and if $p = 2$, then $\tau$ is isotopic to the identity. In general, $\l$ does not admit a diffeomorphism that exchanges $V_1$ and $V_2$ except when $q^2 \equiv \pm 1$ (mod $p$). If $q^2 \equiv 1$ (mod $p$), there exists a diffeomorphism $\sigma_{+}$ that exchanges the Heegaard tori, namely $\sigma_{+}$ : $(u, v) \in V_1 \mapsto (u, v) \in V_2$. If $q = 1$ or $p-1$ then $\sigma_{+}$ is isotopic to the identity. Similarly, when $q^2 \equiv -1$, $\l$ admits a diffeomorphism $\sigma_{-}$ that exchanges $V_1$ and $V_2$, $\sigma_{-}$ : $(u, v) \in V_1 \mapsto (\bar{u}, v) \in V_2$ and $(u, v) \in V_2 \mapsto (u, \bar{v}) \in V_1$.
 
\begin{proof}[Proof of Theorem ~ \ref{class1}] Let us first consider the case of null-homologous knots, the case of not null-homologous knots follows from the same argument. Let $\k$ and $\kp$ be two isotopic null-homologous oriented knots on $\partial V_2$ in $\l$. We always fix the Heegaard torus $\partial V_2$ in $\l$ up to isotopy from the very beginning by using the Theorem 1 of \cite{Bon} so that the knots are homologous on the Heegaard torus too. Since $\k$ and $\kp$ have isomorphic groups, from Lemma ~ \ref{lemma1}, we know that the candidates for $(a', b')$ are $A$, $-A$, $B$, $-B$, $C$, $-C$, $D$ and $-D$. 

We are now going to identify the diffeomorphisms that send $\k$ to possible $\kp$'s and then we are going to analyze when such diffeomorphisms are isotopic to the identity. Clearly, $\tau$ sends $\k$ to $K_{(-a, -b)}$. Note that $\sigma_{+}$ sends $\k$ on $\partial V_2$ to $\k$ on $\partial V_1$. Then after applying the gluing map $\varphi : \partial V_1 \rightarrow \partial V_2$ with $pq' + qp' = 1$, we get
\begin{displaymath} \left(\begin{array}{ccc} -q&q'\\p&p' \end{array}\right) \left(\begin{array}{ccc} a\\b \end{array}\right) = \left(\begin{array}{ccc} -qa + q'b\\pa +p'b \end{array}\right) = \left(\begin{array}{ccc} a'\\b' \end{array}\right) \end{displaymath} 
Note that for $a' = -qa + q'b$ and $b' = pa + p'b$ we have $pa'+ qb' = p(-qa+q'b) + q(pa + p'b) = (pq'+qp')b = b$. By Proposition ~ \ref{prop1}($2$), it follows that we are in the case when $b = pa' + qb'$ and $|b'| = |pa + qb|$. More precisely, $b = pa' + qb'$ and $b' = pa + p'b = pa + qb$ or $b = pa' + qb'$ and $b' = pa + p'b = -pa - qb$. Since $(a, b) = 1$, the latter case does not occur. If we choose $p'$ such that $qp' \equiv 1$ (mod $p$), then we are left with the only case $b = pa' + qb'$ and $b' = pa + qb$ and in this case $(a', b') = (\frac{b-qpa-q^2b}{p}, pa+qb) = C$. Therefore, $\sigma_{+}$ sends $\k$ to $K_C$. In homology one has $[\k] = [K_C]$ that is $b[\lambda_2] = (pa + qb)[\lambda_2] = qb[\lambda_2]$. So the case $C$ occurs only if $q = 1$. Now, $\tau \circ \sigma_{+}$ sends $\k$ to $K_{-C}$ and by homological reasons the case $-C$ occurs only if $q = p-1$.

\par By a similar argument, the diffeomorphism $\sigma_{-}$ sends $\k$ to $K_{D}$. Moreover, in homology $[\k] = [K_D] = b[\lambda_2] = (-pa - qb)[\lambda_2] = - qb[\lambda_2]$. So the case $D$ occurs only if $q = p - 1$. And $\tau \circ \sigma_{-}$ sends $\k$ to $K_{-D}$ and by similar homological reasons the case $-D$ occurs only if $q = 1$. However, by Theorem ~\ref{thmlens}, the diffeomorphism $\sigma_{-}$ exists only when $q^2 \equiv -1$ (mod $p$), for this reason the knot $\k$ is not isotopic to $K_D$ or $K_{-D}$ via a diffeomorphism which is isotopic to the identity. 

\par We want to remark that there is no diffeomorphism of $\l$ sending $\k$ to $K_B$ or $K_{-B}$. We see that such a diffeomorphism cannot be $\sigma_+$ or $\sigma_-$. It cannot be $\tau$ either since $\tau^2 = id$ gives us a contradiction. 

\par Now, using Theorem ~\ref{thmlens} we observe that when $p = 2$ the diffeomorphisms $\tau$ and $\sigma_{+}$ are isotopic to the identity and hence we have Case $(3)$. In Case $(2)$, when $p \neq 2$ and $q = 1$ or $p-1$ the knots $\k$ and $K_C$ are isotopic since in this case only $\sigma_{+}$ is isotopic to the identity. In the remaining cases, only $\tau$ exists and when $p \neq 2$, $\tau$ is not isotopic to the identity. This proves Case $(1)$.
\end{proof}

\begin{lemma} A torus knot $\k$ in $\l$ has a rational Seifert surface $S_{\k}$ of Euler characteristic
\begin{equation*}
\chi({S_{\k}}) = \frac{|rb| + (1 - |rb|)|rap + rbq|}{p}
\end{equation*}
where $r$ is the order of $\k$.
\label{rationalemma}
\end{lemma}
See {\cite[Lemma 2.2]{BJE}} for rational Seifert surface construction for torus knots.
\begin{proof} Let $\k$ be a rationally null-homologous torus knot of order $r$ in $\l$. Considering the corresponding meridional curves of the Heegaard splitting on $\partial V_2$, for any torus knot $\k$ of order $r$ we have 
\begin{equation*}
r[K] = \mathit{m}[\mu_1] + \mathit{l}[\mu_2]= \mathit{m}(-q[\mu_2] + p[\lambda_2]) + \mathit{l}[\mu_2]
\end{equation*}
where $m = \frac{rb}{p}$, $l = ra + mq = ra + \frac{rb}{p}q$ and $p \mid rb$.

We may construct a rational Seifert surface $S_{\k}$ for $r$ copies of $\k$ by taking $|m|$ parallel copies of the meridional disk $\mu_1$ of $\partial V_1$ and $|l|$ parallel copies of the meridional disk $\mu_2$ of $\partial V_2$ and then attaching a half twisted band at each intersection for a total number of $p|l||m| = |l||rb|$ bands. Then, the Euler characteristic $\chi({S_{\k}})$ of $S_{\k}$ is $\chi({S_{\k}}) = \#(disks) - \#(bands)$:
\begin{align*}
\chi({S_{\k}}) &= |l| + |m| - |l||rb| \nonumber \\ 
&= \smallskip |ra + \frac{rbq}{p}| + \frac{|rb|}{p} - |ra + \frac{rbq}{p}||rb| \nonumber \\ 
&= \medskip \frac{|rb| + (1 - |rb|)|rap + rbq|}{p}. \nonumber 
 \qedhere
\end{align*}
\end{proof}


\section{Convex surfaces}
An oriented smooth surface $\Sigma$ in a contact $3$--manifold is called \textit{convex} if there is a \textit{contact vector field} $v$, that is a vector field whose flow preserves the contact structure $\xi$, transverse to $\Sigma$. Given a convex surface $\Sigma$ in a contact $3$-manifold with a contact vector field $v$, the \textit{dividing set} $\Gamma_{\Sigma}$ of $\Sigma$ is defined as $\Gamma_{\Sigma} =\{ x \in \Sigma : v(x) \in \xi_x\}$. The dividing set $\Gamma_{\Sigma}$ is a multi-curve and possibly disconnected. The dividing set $\Gamma_{\Sigma}$ is transverse to the characteristic foliation, splits $\Sigma$ into two subsurfaces $\Sigma \setminus \Gamma_{\Sigma} = \Sigma_{+} \sqcup \Sigma_{-}$ and there is a vector field $w$ that expands/contracts a volume form $\omega$ on $\Sigma_{+}$/ $\Sigma_{-}$ and $w$ points out of $\Sigma_{+}$.  

\begin{theorem}[Giroux's tightness criterion] A convex surface $\Sigma$ in a contact $3$-manifold has a tight neighborhood if and only if $\Sigma \neq S^2$ and $\Gamma_{\Sigma}$ has no homotopically trivial dividing curves or $\Sigma = S^2$ and $\Gamma_{\Sigma}$ is connected.
\end{theorem}

For more information and details, see \cite{Giru, Ko}.

\subsection{Legendrian knots.} 
\begin{figure}[ht]
\begin{center}
     \includegraphics[width=0.6\textwidth]{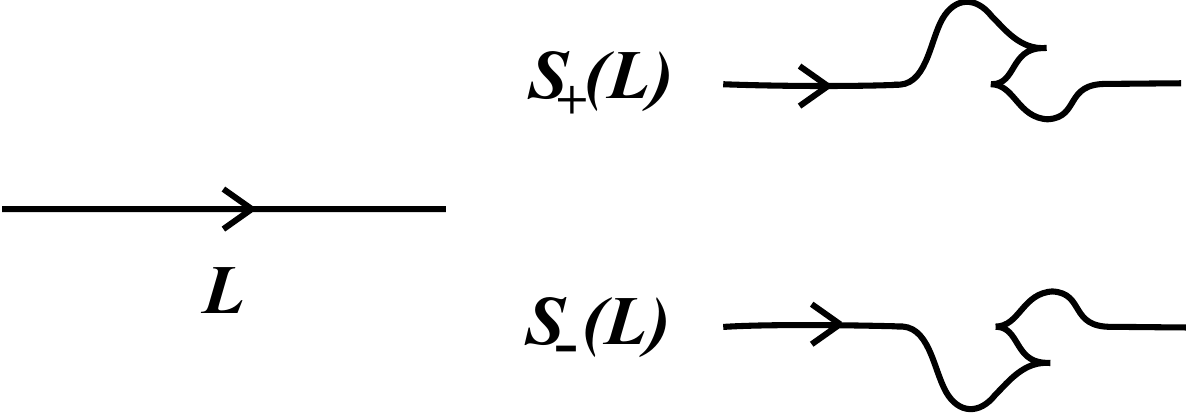}
 \caption{A positive and a negative stabilization of $L$.}
  \label{stab}
\end{center}
\end{figure}

The \emph{positive stabilization/negative stabilization} $S_{+}(L)$/$S_{-}(L)$ of a Legendrian knot $L$ in the standard tight contact structure on $\mathbb{R}^3$ is obtained by modifying the front projection of $L$ by adding a down cusp/an up cusp as in Figure~\ref{stab}, respectively. Since stabilizations are performed locally, by Darboux's theorem this defines stabilizations of Legendrian knots in any contact $3$--manifold. 

\begin{proposition}[{Kanda \cite{Kanda}}] Let $L$ be a Legendrian curve on a surface $\Sigma$ and let $\mbox{\tt tw}_{\Sigma}(L)$ denote the twisting of the contact planes along $L$ measured with respect to the framing on $L$ given by $\Sigma$. Then $\Sigma$ may be made convex relative to $L$ if and only if $\mbox{\tt tw}_{\Sigma}(L) \leq 0$. If $\Sigma$ is a convex surface with dividing curve $\Gamma$, then
\begin{equation}
\mbox{\tt tw}_{\Sigma}(L) = -\frac{1}{2}\# (L \cap \Gamma)
\end{equation}
where $\# (L \cap \Gamma)$ is the unsigned count of intersection number of $L$ and $\Gamma$. Moreover, if $\Sigma$ is a Seifert surface of a single oriented Legendrian curve $L$, the above formula computes the Thurston-Bennequin invariant $\mbox{\tt tb}(L)$ of $L$ and in this case the rotation number $\mbox{\tt rot}(L)$ of $L$ is 
\begin{equation}
\mbox{\tt rot}(L) = \chi(\Sigma_{+}) - \chi(\Sigma_{-}).
\end{equation}
\label{thmKanda}
\end{proposition}

From Proposition ~\ref{thmKanda} we have

\begin{lemma} A surface $\Sigma$ with Legendrian boundary may be made convex if and only if the twisting of contact planes along each boundary component is less than or equal to zero.
\label{lemmaKanda}
\end{lemma}

\subsection{Convex torus in standard form.}  For relatively prime integers $a$, $b$, the \emph{slope} of an $(a,b)$--curve on a torus is $\frac{b}{a}$. A \emph{convex torus (in standard form) with slope $s$} is a torus whose characteristic foliation consists of $2n$ lines of singularities with slope $s$, called \emph{Legendrian divides} and the rest of the foliation is by non-singular lines of slope $r \neq s$, called \emph{Legendrian rulings} where $r$ and $s$ are rational numbers. The $2n$ curves of the dividing set lie between the Legendrian divides. By Giroux's flexibility theorem, \cite{Giru, Ko}, any convex torus with slope $s$ in a tight contact $3$--manifold can be put in a standard form with any ruling slope $r \neq s$.

\begin{theorem}[Classification of tight contact structures on a solid torus, \cite{Ko}] There are $|(r_{0}+1)\cdots (r_{k-1}+1)(r_{k})|$ tight contact structures on a solid torus $S^1 \times D^2$ with standard convex boundary having two dividing curves of slope $-\frac{p}{q}$, where $p > q > 0$ and $-\frac{p}{q} = r_{0} - \frac{1}{r_{1}- \frac{1}{r_{2} \cdots -\frac{1}{r_k}}}$ for $|r_i| < -1$. Moreover, all these contact structures are distinguished by the number of positive regions on a convex meridional disk with Legendrian boundary.
\label{solidtorus}
\end{theorem}

\begin{proposition}[{\cite[Proposition 4.16]{Ko}}] Let $\xi$ be a tight contact structure on $T^2 \times I$ with convex boundary having boundary slopes $s_0$ and $s_1$ on the boundary. Then for any $s$ between $s_0$ and $s_1$, there is a convex torus parallel to the boundary of $T^2 \times I$ with slope $s$.
\label{lemma35}
\end{proposition}

\subsection{Bypasses} 
\begin{figure}[ht]
\begin{center}
  \includegraphics[width=0.6\textwidth]{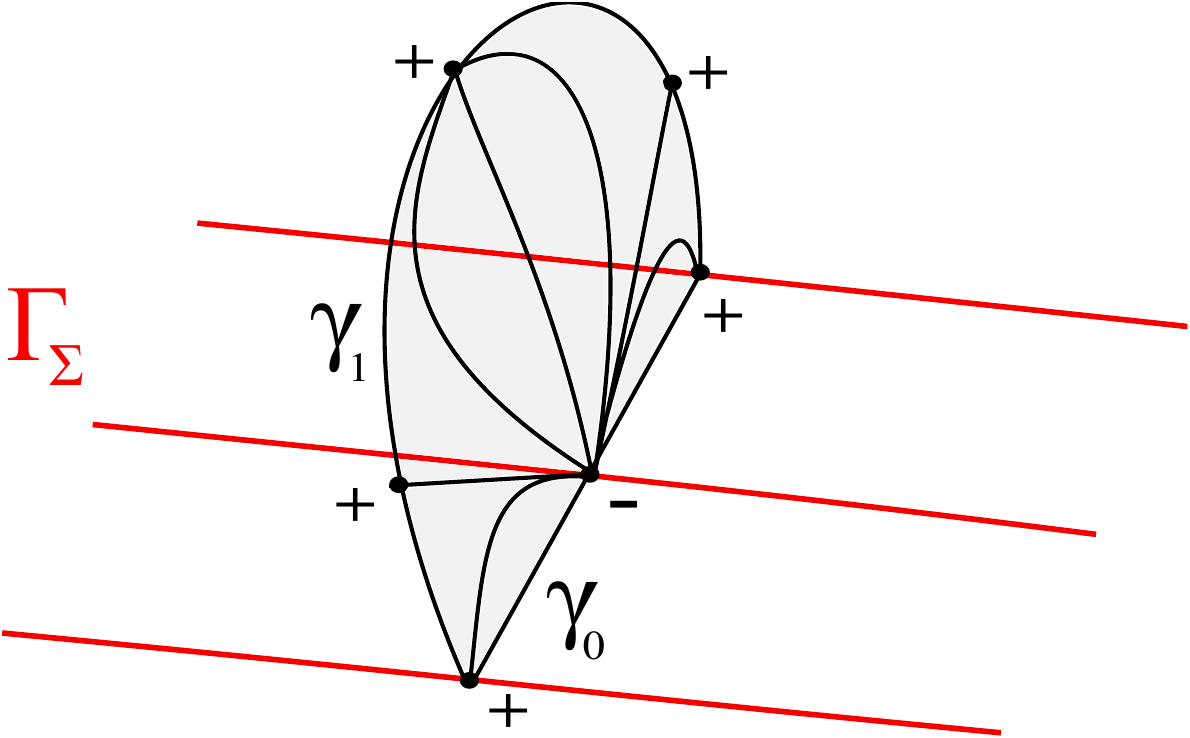}
 \caption{A bypass disk $D$.}
  \label{bypass}
\end{center}
\end{figure}

Let $\Sigma$ be a convex surface in a contact $3$--manifold, a bypass for $\Sigma$ is a convex half disk $D$ (or $D$ with opposite orientation) with Legendrian boundary such that 
\begin{enumerate} 
\item $\partial D = \gamma_{0} \cup \gamma_{1}$ where $\gamma_{0}$, $\gamma_{1}$ are two arcs that intersect at their end points,
\item $D \cap \Sigma = \gamma_{0}$,
\item the characteristic foliation of $D$ has three elliptic singularities along $\gamma_{0}$, two positive elliptic singularities at the end points of $\gamma_{0}$ and one negative elliptic singularity on the interior of $\gamma_{0}$, and only positive singularities along $\gamma_{1}$, alternating between positive elliptic and positive hyperbolic singularities,
\item $\gamma_{0}$ intersects $\Gamma_{\Sigma}$ exactly at the three elliptic singularities of $\gamma_{0}$.
\end{enumerate}
The \emph{sign} of a bypass disk is defined to be the sign of the singularity at the center of the half disk. Figure ~\ref{bypass} is a diagram illustrating a bypass disk.

A dividing curve $\gamma \subset \Gamma_{\Sigma}$ is called \emph{boundary parallel} if $\gamma$ cuts off a half disk of $\Sigma$ which contains no other component of $\Gamma_{\Sigma}$ in its interior. A boundary parallel dividing curve allows us to find bypasses. 

\begin{proposition} [{Honda \cite[Proposition 3.18]{Ko}}] Let $\Sigma$ be a convex surface with Legendrian boundary. If the dividing set $\Gamma_{\Sigma}$ contains a boundary parallel dividing curve $\gamma$, then there exists a bypass for $\Sigma$, provided that $\Sigma$ is not a disk with $\mbox{\tt tb}(\partial\Sigma) = -1$.
\label{bp}
\end{proposition}

\begin{proposition}[{Imbalance Principle, \cite[Proposition 3.17]{Ko}}] Let $\Sigma = S^1 \times [0,1]$ be a convex annulus with Legendrian boundary embedded in a tight contact $3$-manifold. If $tw_{\Sigma}(S^1\times\{0\}) < tw_{\Sigma}(S^1\times\{1\}) \leq 0$, then there exists a boundary parallel curve and hence a bypass along $S^1\times\{0\}$.
\label{ibp}
\end{proposition}


\section{Legendrian torus knots in lens spaces}

A tight contact structure on a $3$--manifold is \emph{universally tight} if its pullback to the universal cover is tight. In this section, we classify oriented Legendrian torus knots $\L$ of knot type $\k$ for $a,b > 0$ in universally tight contact structures on a lens space $\l$. Legendrian torus knots $\L$ with $a,b > 0$ are called \emph{positive} Legendrian torus knots.

There are exactly two tight contact structures on $\l$ with $q \neq p - 1$ which are universally tight, and there is only one if $q = p -1$, {\cite[Proposition 5.1(3)]{Ko}}. According to \cite{Ko}, we can express $\l$ as the union of two solid tori $V_1$ and $V_2$ where $\partial V_1$ has dividing curves of slope $\infty$. Then, we split $V_2$ into a thickened torus $T^2 \times [0, 1]$ with two dividing curves of slope $s_0 = 0$ on $T^2 \times \{0\}$ and $s_1 = - \frac{p}{q}$ on $T^2 \times \{1\}$ and a solid torus with a unique tight contact structure on it. The universally tight contact structures on $\l$ are induced from the universally tight contact structures on $T^2 \times [0, 1]$. By Proposition 5.1 (1) of \cite{Ko} there are two such universally tight contact structures on $T^2 \times [0, 1]$ and they satisfy $PD(e(\xi, s)) = \pm((-q,p)-(-1,0))$. We assume that $\xi_{ut}$ is induced from the one with $PD(e(\xi, s)) = (-q,p)-(-1,0)$. The results in this section similarly hold for the other case and can be easily written down.

\begin{remark} By Proposition ~\ref{lemma35}, in a universally tight lens space $\l$, one can find a convex Heegaard torus $T$ with two dividing curves of any slope in $(-\frac{p}{q}, 0)$. \label{remark1}
\end{remark}

\par Let $\L$ be a Legendrian torus knot of knot type $\k$ of order $r$ in $\xi_{ut}$ on $\l$. We define the rational Legendrian knot invariants which are defined and studied in \cite{F, BG, EB, GO} for rationally null-homologous knots. By Equation $(1)$ in Proposition ~\ref{thmKanda}, by using the set of dividing curves $\Gamma$ for the Heegaard torus containing $\L$, the twisting of $\L$ is $-\frac{1}{2} \# (\L \cap \Gamma)$. Since $\L$ is rationally null-homologous, $\L$ has a framing given by a rational Seifert surface. By using the rational Seifert surface constructed in Lemma ~ \ref{rationalemma}, we compute the rational Seifert framing of $\L$ as $\frac{1}{r}\frac{plm}{r} = \frac{1}{r}lb$ where $m =\frac{rb}{p}$, $l = ra + \frac{rb}{p}q$. The \emph{rational Thurston-Bennequin invariant} of $\L$ is defined as the twisting of $\L$ with respect to the rational Seifert framing and it is denoted by $ \mbox{\tt tb}_{\mathbb{Q}}(\L )$.

Note that an arbitrary $(a,b)$--curve and $(c,d)$--curve on a torus intersect $|det\left(\begin{array}{ccc} a&c\\b&d \end{array}\right)|$ times. If the dividing curves $\Gamma$ of the Heegaard torus containing $\L$ have slope $ -\frac{s}{t}$ for integers $s,t > 0$ and if $2n$ is the number of dividing curves, then the \emph{rational Thurston-Bennequin invariant} of $\L$ is 
\begin{equation*}
\mbox{\tt tb}_{\mathbb{Q}}(\L) = \frac{1}{r}lb - n|as + bt|.
\end{equation*}

\par Let $\mathcal{L}(K)$ denote the set of all rationally null-homologous Legendrian knots in knot type $K$. The \emph{maximal rational Thurston-Bennequin invariant} $\mbox{\tt tb}_{\mathbb{Q}}(K)$ of the knot type $K$ is defined as 
\begin{equation*}
\mbox{\tt tb}_{\mathbb{Q}}(K) =  max\{\mbox{\tt tb}_{\mathbb{Q}}(L) \mid L \in  \mathcal{L}(K)\}.
\end{equation*}

\begin{theorem} For $a, b > 0$ relatively prime integers, the maximal rational Thurston-Bennequin invariant $\mbox{\tt tb}_{\mathbb{Q}}(\k)$ is 
\begin{equation*}
\mbox{\tt tb}_{\mathbb{Q}}(\k) = ab - a - b + b^2 q/p.
\end{equation*}
\label{tbmax}
\end{theorem}
\begin{proof} By Remark ~\ref{remark1}, we can find a convex Heegaard torus $T$ with two dividing curves of any slope in $(-\frac{p}{q}, 0)$. In particular, there is a convex Heegaard torus $T$ with two dividing curves of slope $-1$. For $a, b > 0$ relatively prime integers, isotope the Legendrian ruling curves on $T$ to have slope $\frac{b}{a}$ so that the ruling curves are Legendrian torus knots $\L$ of knot type $\k$. Since the intersection number $\# (\L \cap \Gamma)$ is minimal on this convex torus $T$, the rational Thurston-Bennequin invariant is maximal. The maximal rational Thurston-Bennequin invariant of the knot type $\k$ computed as
 $$\mbox{\tt tb}_{\mathbb{Q}}(\k) = \frac{1}{r}lb - |a + b|$$
where $l = ra + \frac{rb}{p}q$. 
\end{proof}
\par The \emph{rational rotation number} $\mbox{\tt rot}_{\mathbb{Q}}(L)$ of an oriented rationally null-homologous Legendrian knot $L$ of order $r$ is defined as the winding number of $TL$ after trivializing the contact structure along a rational Seifert surface for $L$ divided by $r$. 

\par Let $\L$ be an oriented Legendrian torus knot of order $r$ with maximal rational Thurston-Bennequin invariant that sits on a convex Heegaard torus $T$. In what follows, we will explain how to compute the rational rotation number of $\L$ in a similar way as Etnyre and Honda computed for Legendrian torus knots in the standard tight $S^3$, \cite{EH}.

\par Let $L(p, q) = V_1 \cup _{T} V_2$ where $V_1$ and $V_2$ are both $D^2\times S^1$ with meridional curve $\mu_1$ and $\mu_2$ respectively. Define an invariant $f_T$ of homology classes of curves on a convex Heegaard torus $T$ as follows: Let $v$ be any globally nowhere zero section of $\xi_{ut}$ and $w$ a nowhere zero section of $\xi_{ut}|_T$ which is tangent to the Legendrian divides and transverse to and twists with $\xi$ along the Legendrian ruling curves. Let $f_T(\gamma)$ equal the rotation of $v$ relative to $w$ along a closed oriented curve $\gamma$ on $T$. For details and the properties of the function $f_T$, see \cite{Et} and \cite{EH}. If $L$ is a Legendrian ruling or a Legendrian divide of order $r$ on $T$ then $f_T(L) = r\, \mbox{\tt rot}_{\mathbb{Q}}(L)$. The rational rotation number of $\L$ of order $r$ on the Heegaard torus $\partial V_2 = T$ is 
$$r\,\mbox{\tt rot}_{\mathbb{Q}}(\L) = m f_T(\mu_1) + l f_T(\mu_2) $$
where $m = \frac{rb}{p}$, $l = ra + \frac{rb}{p}q$. 


\begin{theorem} Let $\L$ be an oriented Legendrian torus knot of order $r$ with maximal rational Thurston-Bennequin invariant. If $a, b > 0$, then the range of possible rational rotation numbers $\mbox{\tt rot}_{\mathbb{Q}}(\L)$ is 
$$\{\pm b(1 - \frac{q+1}{p})\}.$$
\label{postorusk}
\end{theorem}
\begin{proof} By the proof of Theorem ~\ref{tbmax}, a Legendrian knot $\L$ with maximal rational Thurston-Bennequin invariant is on a convex Heegaard torus $T$ with two dividing curves of slope $-1$ in $\l = V_1 \cup_T V_2$. To compute the rational rotation number for $\L$, we need to compute $f_T(\mu_1)$ and $f_T(\mu_2)$. 

$f_T(\mu_1) = \pm(p-q-1)$: Recall that the meridional curve $\mu_1 = \partial{D_{V_1}}$ is a $(-q, p)$--curve on $T = \partial{V_2}$. Isotope the Legendrian ruling curves on $T$ to be $(-q, p)$--curves so that $\mu_1$ is a ruling curve. In this case, the twisting of the contact planes along $\mu_1$ is $-(p-q) < 0$. So, by Lemma ~\ref{lemmaKanda} we may isotope $D_{V_1}$ to be convex. By the proof of Theorem ~\ref{tbmax}, the dividing curves on $T$ have slope $-1$ and hence they are $(-1,1)$--curves and intersect $\mu_1$ $2(p-q)$ times. So, the dividing curves on $D_{V_1}$ intersect $\mu_1$ $2(p-q)$ times. By following the proof of Fact 1 of \cite{EH}, the dividing curves on $D_{V_1}$ separate off disks of the same sign that contain no other dividing curves. Then, by Equation $(2)$ of Proposition ~\ref{thmKanda}, we have $f_T(\mu_1) = (p-q)-1$ or $f_T(\mu_1) = 1-(p-q)$.

$f_T(\mu_2) = 0$: Isotope the Legendrian ruling curves on $T = \partial{V_2}$ to be meridional. Let $D_{V_2}$ be the meridional disk of $V_2$ where $\mu_2 = \partial{D_{V_2}}$ is a ruling curve. Since the twisting of the contact planes along $\mu_2$ is $-1$, by Lemma ~\ref{lemmaKanda}, we can isotope $D_{V_2}$ to be convex. We know that the dividing curves on $D_{V_2}$ intersect $\mu_2$ twice. Thus, we have only one possible configuration for the dividing curves of $D_{V_2}$. By Equation $(2)$ of Proposition ~\ref{thmKanda}, $f_T(\mu_2) = 0$. 
\end{proof}

\begin{theorem} \label{thmain} Two oriented positive Legendrian torus knots in a universally tight contact structure on a lens space are uniquely realized up to contactomorphism if and only if their oriented knot types, rational Thurston-Bennequin invariants and rational rotation numbers agree.
\end{theorem} 

The following two lemmas provide us the proof of the Theorem ~\ref{thmain}.

\begin{lemma} Two oriented Legendrian $(a,b)$--torus knots, $a,b >0$, $L$ and $L'$ in $(\l, \xi_{ut})$ with maximal rational Thurston-Bennequin invariant are uniquely realized up to contactomorphism if and only if $\mbox{\tt rot}_{\mathbb{Q}}(L) = \mbox{\tt rot}_{\mathbb{Q}}(L')$.
\end{lemma}
\begin{proof} Let $T$ and $T'$ be standard convex Heegaard tori on which $L$ and $L'$ respectively sit in $L(p, q)$. Also, let $V_1 \cup _{T} V_2$ and $V_1' \cup _{T'} V_2'$ be the Heegaard splittings associated to $T$ and $T'$. Since $\mbox{\tt tb}_{\mathbb{Q}}(L) = \mbox{\tt tb}_{\mathbb{Q}}(L') = \mbox{\tt tb}_{\mathbb{Q}}(\k)$, the slopes of the dividing curves on $T$ and $T'$ are the same. Then, by Theorem ~\ref{solidtorus}, by the classification of tight contact structures on solid tori, there is a contactomorphism $g : V_1 \rightarrow V_1'$ such that $g(L) = L'$. By Theorem ~\ref{solidtorus} again, the contactomorphism type of a tight contact structure on $V_2$ or $V'_2$ is determined by the number of positive bypasses on meridional disks. If $r$ is the order of $L$ and $L'$ in $\l$, then the number of positive bypasses on meridional disks are determined by $r$ times the rational rotation number of the Legendrian knots $L$ and $L'$, respectively. We can extend the contactomorphism $g$ to all of $L(p, q)$ provided that $L$ and $L'$ have the same rational rotation number.
\end{proof}
\begin{lemma} If $\L$ is a positive Legendrian torus knot in $(\l, \xi_{ut})$ with non-maximal rational Thurston-Bennequin invariant then there is a Legendrian torus knot $\Lp$ such that $\L$ is a stabilization of $\Lp$.
\end{lemma}
\begin{proof} Let $T$ be a standard convex Heegaard torus on which the positve Legendrian torus knot $\L$ sits. Let $-\frac{s}{t}$ for $s, t > 0$ be the slope of the the dividing curves $\Gamma_T$ on $T$ and let $2n$ be the number of dividing curves. By Theorem ~\ref{tbmax}, Legendrian torus knots with maximal rational Thurston-Bennequin invariant sit on a convex Heegaard torus with two dividing curves of slope $-1$. Since $\mbox{\tt tb}_{\mathbb{Q}}(\L) < \mbox{\tt tb}_{\mathbb{Q}}(\k)$, we have two cases for the slope and the number of dividing curves of $\Gamma_T$ on $T$: $-\frac{s}{t} = -1$ and  $n > 1$ or $-\frac{s}{t} \neq -1$ and $n \geq 1$.

We are now going to show that we can find a bypass disk in both cases. By Remark ~\ref{remark1}, we know that there is a convex torus $T'$ with two dividing curves of slope $-1$. In fact, by \cite{Ko} we can assume that $T'$ is a standardly embedded convex torus parallel and disjoint from $T$. Now take the $T^2 \times [0,1]$ region between $T$ and $T'$ and take the annulus $A = \L \times [0,1]$ between $T = T^2 \times \{0\}$ and $T' = T^2 \times \{1\}$. Furthermore, isotope the ruling curves on both boundary components of $T^2 \times [0,1]$ to have slope $\frac{b}{a}$. Then $\partial A = \L \cup \Lp $ are Legendrian ruling curves on the boundary of $T^2 \times [0,1]$ and the twisting of contact planes along each boundary component will be less than zero. Therefore, by Lemma ~\ref{lemmaKanda}, we can make $A$ convex. The dividing curves on $T = T^2 \times \{0\}$ are $(-t,s)$--curves and intersect $A$ in $2n|det\left(\begin{array}{ccc} a&-t\\b&\phantom{-}s \end{array}\right)|=2n(sa+tb)$ points and the dividing curves on $T' = T^2 \times \{1\}$ are $(-1,1)$--curves and intersect $A$ in $2|det\left(\begin{array}{ccc} a&-1\\b&\phantom{-}1 \end{array}\right)| = 2(a+b)$ points. In both cases when $-\frac{s}{t} = -1$ and $n > 1$ or when $-\frac{s}{t} \neq -1$ and $n \geq 1$, we have $2n(sa+tb) > 2(a+b)$. So, there is a boundary parallel dividing curve along $T = T^2\times \{0\}$ and hence by Proposition~\ref{ibp} a bypass for $\L$. In other words, $\L$ destabilizes.
\end{proof}

\end{document}